\documentclass{article}

\usepackage{amssymb}
\usepackage{amsmath,ifthen,amsthm,cite}

\newcommand{\HH}{\mathcal{H}}
\newcommand{\C}{\mathbf{C}}
\newcommand{\N}{\mathbf{N}}
\newcommand{\M}{\mathcal{M}}

\newcommand{\Poincare}{Poincar\'e}

\DeclareMathOperator{\SL}{SL}
\newcommand{\Z}{\mathbf{Z}}
\newtheorem{theorem}{Theorem}

\newtheorem{definition}[theorem]{Definition}

\newtheorem{question}[theorem]{Question}

\newcommand{\sltz}{\SL_2(\Z)}
\newcommand{\smatrix}[4]{\left(\begin{smallmatrix}#1 & #2 \\ #3 & #4\end{smallmatrix}\right)}

\newcommand{\abcd}{\left(\begin{matrix}a & b \\ c & d\end{matrix}\right)}
\begin{document}

\title{Generating spaces of modular forms with~$\eta$-quotients}
\author{L. J. P. Kilford}

\maketitle

\begin{abstract}
In this note we consider a question of Ono, concerning which spaces of classical modular forms can be generated by sums of $\eta$-quotients. We give some new examples of spaces of modular forms which can be generated as sums of $\eta$-quotients, and show that we can write all modular forms of level~$\Gamma_0(N)$ as rational functions of $\eta$-products.
\end{abstract}

\section{Introduction}

Let~$z$ be an element of the \Poincare{} upper half plane~$\HH$, and let~$q:=\exp(2 \pi i z)$. We recall the definition of the Dedekind~$\eta$-function:
\[
\eta(q):=q^{1/24}\prod_{n=1}^\infty (1-q^n).
\]
It is well-known that this function plays an important role in the theory of modular forms; for instance, the unique normalised cusp form of level~1 and weight~12, the~$\Delta$ function, can be written as~$\Delta(q)=\eta(q)^{24}$. 

Another classical application of the~$\eta$-function is to the theory of \emph{partitions}; a partition of~$n$ is a way of writing~$n$ as a sum of positive integers; for instance,~$4=2+2$ is a partition of~4. The number of partitions of~$n$, $p(n)$, grows quickly with~$n$. For instance, in the early years of the $20^{\rm th}$ century, MacMahon computed that~$p(200)=397299029388$. The reciprocal of the $\eta$-function (with the $q^{1/24}$ removed) gives the following generating function for the~$p(n)$:
\[
\sum_{n=0}^\infty p(n) q^n = \prod_{m=1}^\infty \left(\frac{1}{1-q^m}\right).
\]

It can be shown that~$\eta(q)^{24}$ is a modular form for~$\sltz$ In~\cite{apostol}, an explicit transformation formula for~$\eta(q)$ under the action of~$\sltz$ is given; it is of the form
\begin{eqnarray}
\label{eta-transformation-formula}
\eta\left(\frac{az+b}{cz+d}\right)=\varepsilon\cdot(cz+d)^{1/2}\eta(z)\text{ for }\abcd \in \sltz,
\end{eqnarray}
where~$\varepsilon^{24}=1$. (The explicit definition of~$\varepsilon$, which depends on~$a,b,c,d$, is complicated, but is given in Theorem~3.4 of~\cite{apostol}).

We see from~\eqref{eta-transformation-formula} that $\eta(q)$ satisfies a ``weight~$1/2$'' transformation formula, which is consistent with the statement that~$\Delta$ is a weight~12 modular form.
\begin{definition}
Let~$N$ be a positive integer, let~$\{r_\delta\}$ be a set of integers, and let~$f$ be a meromorphic function from the Poincar\`e upper half plane to~$\C$ of the form
\[
f(z)=\prod_{0 < \delta | N} \eta\left(q^\delta\right)^{r_\delta}, \text{{\rm where }}q:=\exp(2\pi i z).
\]
We call~$f$ an \emph{$\eta$-quotient}, and if all of the~$r_\delta$ are non-negative we say that~$f$ is an \emph{$\eta$-product}.
\end{definition}
We see that~$\Delta$ is an $\eta$-product; we will see some examples of $\eta$-quotients which are modular forms below.

There are some general theorems which tell us when an $\eta$-quotient is a classical modular form, by telling us how it transforms under elements of certain congruence subgroups.
\begin{theorem}[Ligozat~\cite{ligozat}, quoted in~\cite{mcmurdy-thesis}]
Suppose that~$f(z)$ is an $\eta$-quotient which satisfies the following properties:
\begin{equation}
\label{sum-of-r-delta}
\sum_{0 < \delta | N} \delta \cdot r_\delta \equiv 0 \mod 24
\end{equation}
and
\begin{equation}
\sum_{0 < \delta | N} \frac{N}{\delta} \cdot r_\delta \equiv 0 \mod 24.
\end{equation}
Then~$f(z)$ satisfies
\begin{eqnarray}
\label{how-eta-quotients-transform}
f\left(\frac{az+b}{cz+d}\right) = \chi(d) (cz+d)^k f(z)
\end{eqnarray}
for every~$(\begin{smallmatrix}a&b\\c&d\end{smallmatrix}) \in \Gamma_0(N)$, where~$k := \frac{1}{2}\sum_{0 < \delta | N} r_\delta$, and
\begin{equation}
\label{character-formula}
\chi(d):=\left(\frac{(-1)^k\cdot s}{d}\right), \text{{\rm where }} s:=\prod_{0 < \delta | N} \delta^{r_\delta}.
\end{equation}
\end{theorem} 

In~\cite{ono-web-of-modularity}, the following question is posed:
\begin{question}[Ono~\cite{ono-web-of-modularity}, Problem~1.68]
Which spaces of classical modular forms can be generated by sums of $\eta$-quotients?
\end{question}
It is well-known (see for instance~\cite{serre-cours}, Corollary~2 to Theorem~VII.4) that every modular form of level~1 is a polynomial in the Eisenstein series~$E_4$ and~$E_6$, so if we can write~$E_4$ and~$E_6$ in terms of $\eta$-quotients then we can write every modular form of level~1 in terms of $\eta$-quotients. In~\cite{ono-web-of-modularity}, we see that
\[
E_4(z)=\frac{\eta(q)^{16}}{\eta(q^2)^8}+2^8\cdot\frac{\eta(q^2)^{16}}{\eta(q)^8}
\]
and
\[
E_6(z)=\frac{\eta(q)^{24}}{\eta(q^2)^{12}}-2^5\cdot 3 \cdot 5 \cdot \eta(q^2)^{12} - 2^9\cdot 3 \cdot 11 \cdot \frac{\eta(q^2)^{12}\eta(q^4)^8}{\eta(q)^8}+2^{13}\cdot \frac{\eta(q^4)^{24}}{\eta(q^2)^{12}}.
\]
One can check that these are identities of modular forms by verifying that each of these $\eta$-quotients is a modular form on~$\Gamma_0(4)$ and that the $q$-expansions of both sides agree up to a suitable bound, such as the Sturm bound~\cite{stein-sturm} (also known as the Hecke bound); if a modular form has a  zero at the cusp~$\infty$  of large enough degree, then it must be the zero modular form. In this case, we find that the dimension of the space~$M_4(\Gamma_0(4))$ is~3, and the Sturm bound here is~2. We can use a computer algebra package to check that these really are equalities.

In this paper we will give further examples of this form, and prove a theorem which says that we can write \emph{every} modular form for a congruence subgroup~$\Gamma_0(N)$ in terms of $\eta$-quotients of level at most~$\Gamma_0(4N)$.

\section{Writing eigenforms as $\eta$-products or $\eta$-quotients}

We will briefly consider $\eta$-products, and show that there are only finitely many spaces of modular forms which can be generated by $\eta$-products. These do exist: for instance, $S_2(\Gamma_0(11)) = \C\cdot \eta(q)^2\eta(q^{11})^2$, and~$S_1(\Gamma_1(23)) = \C\cdot \eta(q)\eta(q^{23})$. However, it is clear that only finitely many spaces of modular forms can be completely generated by $\eta$-products, as there are only finitely many $\eta$-products with $q$-expansion of the form~$q^a+O(q^{a+1})$ for any~$a$; this follows from equation~(\ref{sum-of-r-delta}), because as the level increases the sum increases also.

Questions of this sort have been asked before; in~\cite{biagioli}, the question of when an $\eta$ product is also a simultaneous eigenform for the Hecke operators~$T_p$ is addressed. In~\cite{dummit-et-al}, a complete list of the 30 different $\eta$-products which are classical modular cuspidal eigenforms of integer weight with no zeroes outside the cusps is given (two of these have half-integral weight). Another proof of this classification was given by~\cite{koike-eta-products}.

Similarly, in~\cite{martin-eta-quotients}, a complete list of all $\eta$-quotients which are holomorphic new eigenforms of integral weight for the congruence subgroups~$\Gamma_1(N)$ is derived. There are 74 of these, with levels ranging from~$\sltz$ (the $\Delta$-function) up to $\Gamma_1(576)$ (the modular form $\eta(q^{12})^{-2}\cdot\eta(q^{24})^6\cdot \eta(q^{48})^{-2} \in S_1(\Gamma_1(576))$).
There have been other investigations; in~\cite{koehler}, a list of the 65 $\eta$-products which have weight~1 and level~$\Gamma^*(36)$ is given; the congruence subgroup~$\Gamma^*(N)$ is generated by~$\Gamma_0(36)$ and the Fricke involution~$\smatrix{0}{1/\sqrt{N}}{-1/\sqrt{N}}{0}$. Another angle of inquiry about $\eta$-products is given in~\cite{gordon-robins}, which considers $\eta$-products which are lacunary modular forms; in other words, they have many zero Fourier coefficients.
\section{Sums of $\eta$-products and $\eta$-quotients}

We see that, from the results we have mentioned in the previous section, if we are to generate all spaces of modular forms for a given congruence subgroup~$\Gamma$ we will need to expand from $\eta$-products and $\eta$-quotients to consider more general objects. We will consider sums of $\eta$-products and $\eta$-quotients in this section.

From equation~\eqref{how-eta-quotients-transform}, we see that the only characters we can have are the Legendre symbols, so if we are considering congruence subgroups containing~$\Gamma_1(N)$ we are restricted to the cases of~$\Gamma_0(N)$ and~$\Gamma_0(N)$ plus a character of order~2.

We will be using the standard formulae for the dimension of spaces of modular forms in this section; these can be found in many books on modular forms; one reference is Theorem~3.5.1 and Figure~3.3 of~\cite{diamond-shurman}. Let~$g=\lfloor(p+1)/12\rfloor$ if~$p \not\equiv 1 \mod 12$ and~$\lfloor(p+1)/12\rfloor -1$ otherwise (this is the genus of~$X_0(p)$). The formulae (if~$N=p$ is prime) are (for~$k \ge 4$ and~$N \ge 5$):
\begin{eqnarray*}
\dim S_k(\Gamma_0(p)) 	&=& (k-1)\cdot(g-1)+\left\lfloor\frac{k}{4}\right\rfloor\cdot \left(1+\left(\frac{-1}{p}\right)\right)\\
&+&\left\lfloor\frac{k}{3}\right\rfloor\cdot\left(1+\left(\frac{-3}{p}\right)\right)+k-2.
\end{eqnarray*}
The dimension of~$S_2(\Gamma_0(p))$ is~$g$, and the dimension of the space of Eisenstein series is~1 if~$k=2$ and~2 if~$k \ge 4$ (the Eisenstein series of weight~$k \ge 4$ are the images of~$E_k$ under the two degeneracy maps from level~1 to level~p, and the Eisenstein series of weight~2 is defined in terms of the Eisenstein series of weight~2 which is not a modular form). There are similar formulae for~$N$ composite and for other congruence subgroups.

We will first consider some more examples of spaces of modular forms which can be generated by $\eta$-quotients: for example, we can write every modular form of level~$\Gamma_0(2)$ as a sum of $\eta$-quotients. It will suffice to show that~$M_2(\Gamma_0(2))$, $M_4(\Gamma_0(2))$, $M_6(\Gamma_0(2))$ and $M_8(\Gamma_0(2))$ can be generated by sums of $\eta$-quotients, because the following relation holds:
\[
S_{k+8}(\Gamma_0(2)) = f \cdot M_k(\Gamma_0(2)), \text{ for }k \in \N,
\] 
where~$f= \lambda \cdot (\eta(z)\cdot\eta(2z))^8\in S_8(\Gamma_0(2))$. This can be derived from the dimension formulae that we have given above.


Because we can write level~1 modular forms as sums of $\eta$-quotients, and we know (for instance, from Proposition~III.3.19 of~\cite{koblitz}) that the unique normalised cusp-form~$f$ of weight~$8$ and level~$\Gamma_0(2)$ is an $\eta$-quotient, we need only to consider~$M_2(\Gamma_0(2))$, which is 1-dimensional. Using {\sc Magma}~\cite{magma}, and the fact that the Hecke bound for the space~$M_2(\Gamma_0(8))$ is~3, we see that
\[
\frac{\eta(q^2)^{20}}{\eta(q)^8\cdot \eta(q^4)^8} + 16 \cdot \frac{\eta(q^8)}{\eta(q^2)^4} = E_{2,2} \in M_2(\Gamma_0(2)),
\]
so every modular form of level~$\Gamma_0(2)$ can be written as a sum of $\eta$-quotients. We need to consider the space of level~$\Gamma_0(8)$ because the $\eta$-quotients have level~$\Gamma_0(8)$.

Similarly, to show that every modular form for~$\Gamma_0(3)$ is a sum of $\eta$-quotients, we need to check that~$E_{3,2} \in M_2(\Gamma_0(3))$ is a sum of $\eta$-quotients. After some calculation, we find (using the fact that the Hecke bound for~$M_2(\Gamma_0(27))$ is 6) that
\[
36\frac{\eta(q^9)^6}{\eta(q^3)^2} + 12\frac{\eta(q^3)^2}{\eta(q^9)^2}+\frac{\eta(q^9)^{10}}{\eta(q^3)^3\eta(q^{27})^3}+108\frac{\eta(q^9)^2\eta(q^{27})^3}{\eta(q^3)}+9\frac{\eta(q^3)^3\eta(q^{27})^3}{\eta(q^9)^2}+9\frac{\eta(q^{27})^6}{\eta(q^9)^2}
\]
is equal to~$E_{3,2}$; this, together with the knowledge that~$S_6(\Gamma_0(3))=\C \cdot \eta(q)^6\eta(q^3)^6$ and that the Eisenstein series of weights~4 and~6 and level~$\Gamma_0(3)$ are oldforms and therefore the image of forms of level~1 which can be written as $\eta$ quotients, shows us that we can write modular forms of level~$\Gamma_0(3)$ as sums of $\eta$ quotients of level at most~$4\cdot27$.

In fact, we can write~$E_4$ and~$E_6$ as sums of eta-quotients of level~27, so we can actually write every modular form for~$\Gamma_0(3)$ in terms of $\eta$ quotients of level~27.

In a similar fashion, we can write the elements of~$M_3(\Gamma_1(3))$ as a linear combination of~$\eta(q^3)^9/\eta(q)^3$ and~$\eta(q)^9/\eta(q^3)^3$; these, together with the Eisenstein series of level~3 and weights~2,4 and~6 and the cusp form of level~3 and weight~6 generate all of the modular forms of level~3 and odd weight, so we can write every modular form of level~$\Gamma_1(3)$ as a sum of $\eta$-quotients.

We can prove similar results for level~$\Gamma_1(4)$ as well. Using Exercise~III.3.17 of Koblitz's book on modular forms~\cite{koblitz}, we see that every $f \in M_k(\Gamma_1(4))$ can be written as a polynomial in the modular forms~$F:=\eta(q^4)^8/\eta(q^2)^4$, which has weight~2, and~$\theta^2:=\eta(q^2)^{10}/(\eta(q)^4\eta(q^4)^4)$, which has weight~1. As both of these are $\eta$-quotients, we are done; every modular form of level~$\Gamma_1(4)$ can be written as a sum of $\eta$-quotients of level at most~4.

A rather tedious {\sc Magma} calculation shows us that the unique Eisenstein series~$E_{5,2}$ of level~$\Gamma_0(5)$ and weight~2 can be written as the following sum of $\eta$-quotients of level~$\Gamma_0(20)$:
\begin{eqnarray*}
E_{5,2}(q)&=&5\cdot\frac{\eta(q^2)^3\eta(q^4)^4\eta(q^5)^3\eta(q^{10})^5}{\eta(q)^7\eta(q^{20})^4}-2\cdot\frac{\eta(q^2)^7\eta(q^4)^2\eta(q^5)^3\eta(q^{10})}{\eta(q)^7\eta(q^{20})^2} \\
&-&4\cdot\frac{\eta(q^2)^7\eta(q^4)^7\eta(q^{10})^5}{\eta(q)^7\eta(q^5)^5\eta(q^{20})^3}+\frac{\eta(q^4)^4\eta(q^5)^6\eta(q^{10})^4}{\eta(q)^6\eta(q^{20})^4}\\
&+&2\cdot\frac{\eta(q^2)^4\eta(q^4)^2\eta(q^5)^6}{\eta(q)^6\eta(q^{20})^2}-\frac{\eta(q^2)^4\eta(q^4)^7\eta(q^{10})^4}{\eta(q)^6\eta(q^5)^2\eta(q^{20})^3}.
\end{eqnarray*}
This follows because these six $\eta$-quotients generate~$M_2(\Gamma_0(20))$, where~$E_{5,2}$ appears as an oldform. We can find the coefficients of these using linear algebra.

In a very similar way, we can write both elements of~$M_2(\Gamma_0(5),\left(\frac{\cdot}{5}\right))$ as $\eta$-quotients by considering them as oldforms of level~20 and finding a basis of $\eta$-quotients for this space of level~20 forms. This suffices to show that all modular forms with character~$\left(\frac{\cdot}{5}\right)$ can be written as sums of $\eta$-quotients.

At this point, we notice that if~$p \ge 5$ is a prime then we cannot write every modular form for~$\Gamma_0(p)$ as a sum of $\eta$-quotients of level~$p^r$. We use the dimension formulae for~$M_2(\Gamma_0(p))$ to show that there is at least one Eisenstein series of weight~2 with nonzero constant term in its Fourier expansion which we need to write as a sum of $\eta$ quotients, each of which has the form~$\eta(q)^a\eta(q^p)^b\cdots\eta\left(q^{p^r}\right)^z$. Therefore, we need to solve the following two simultaneous equations:
\begin{eqnarray}
a+p\cdot b+p^2\cdot c + \cdots + p^r \cdot z&=&0\\
a+b+c+\cdots +z=4;
\end{eqnarray}
putting these together, we find that
\begin{eqnarray}
\label{reduce-eta-quotient}
(p-1)\cdot b+(p^2-1)\cdot c + \cdots + (p^r-1)\cdot z = -4.
\end{eqnarray}
We now reduce~\eqref{reduce-eta-quotient} modulo~$p-1$; the left hand side vanishes, and the only way that the right hand side can vanish modulo~$p-1$ is if~$p \in \{2,3,5\}$.

Now we consider the case when~$p=5$. We use the formula for the character of the $\eta$-quotient given in~\eqref{character-formula} to work on~\eqref{reduce-eta-quotient}; after dividing both sides by~4, we see that the right hand side is~$-1$, whereas the left hand side has only even entries (because the requirement that the character be trivial forces the sum~$b+d+f+\cdots$ to be even). This is a contradiction.

We can tackle~$\Gamma_0(7)$ in a similar fashion; in this case, we have to check that we can generate~$M_2(\Gamma_0(7))$ and~$M_4(\Gamma_0(7))$ by $\eta$-quotients, and again we prove this by finding a basis of $\eta$-quotients for~$M_2(\Gamma_0(28))$ and~$M_4(\Gamma_0(28))$ and then considering the level~7 spaces as oldforms in these spaces. Here we find that
\begin{eqnarray*}
E_{7,2}&=&16\cdot\frac{\eta(q^2)^6\eta(q^4)^6\eta(q^{14})^2}{\eta(q)^8\eta(q^{28})^2}-32\cdot\frac{\eta(q^2)^4\eta(q^4)^3\eta(q^7)^6}{\eta(q)^6\eta(q^{14})^2\eta(q^{28})}\\
&+&\frac{49}{2}\cdot\frac{\eta(q^4)\eta(q^7)^8\eta(q^{14})^2}{\eta(q)^4\eta(q^{28})^3}-10\cdot\frac{\eta(q^2)^4\eta(q^4)^6\eta(q^{14})^4}{\eta(q)^4\eta(q^7)^4\eta(q^{28})^2}\\
&-&\frac{16}{7}\cdot\frac{\eta(q^2)^5\eta(q^7)^5}{\eta(q)^3\eta(q^{14})}+6\cdot\frac{\eta(q^2)^7\eta(q^4)^2\eta(q^{28})^2}{\eta(q)\eta(q^7)\eta(q^{14})^5}\\
&-&\frac{17}{14}\cdot\frac{\eta(q)\eta(q^2)^3\eta(q^7)^5\eta(q^{14})^5}{\eta(q^4)^7\eta(q^{28})^3};
\end{eqnarray*}
the equations for the three cuspforms in~$S_4(\Gamma_0(7))$ are similar.


\begin{table}
\begin{center}
\begin{tabular}{|c|c|c|}
\hline
Group & Level of $\eta$-products needed\\
\hline
$\SL_2(\Z)$ &  4 or 27\\
\hline
$\Gamma_1(2)$ & 8\\
\hline
$\Gamma_1(3)$ & $27$\\
\hline
$\Gamma_1(4)$ & 4\\
\hline
$\Gamma_0(5)$ & 20\\
\hline
$\Gamma_0(7)$ & 28\\
\hline
\end{tabular}
\caption{Levels of $\eta$-quotients needed to generate spaces of modular forms}
\label{table-of-spaces}
\end{center}
\end{table}

\section{Generating spaces of modular forms with rational functions of $\eta$-quotients}

We now consider a larger set of functions than in the previous section; we will now consider \emph{rational functions} of $\eta$-quotients and $\eta$-products. 
We will now prove the following theorem on the generation of spaces of modular forms by rational functions of $\eta$-quotients (and in fact by rational functions of $\eta$-products):
\begin{theorem}
Let~$N$ be a positive integer, let~$k$ be a positive even integer, and let~$f \in M_k(\Gamma_0(N))$ be a modular form. Then we can write~$f$ as a rational function of~$\eta$-quotients.
\end{theorem}
\begin{proof}
We first note that the classical $j$-invariant is an $\eta$-quotient, because we can write it as
\[
j:=\frac{E_4^3}{\Delta},
\]
where~$E_4$ is an $\eta$-quotient and~$\Delta=\eta(q)^{24}$. This means that the modular function~$j(q^N)$ is also an $\eta$-quotient, for any~$N \in \N$.

We now quote a theorem about the field of modular functions of weight~0 for the congruence subgroup~$\Gamma_0(N)$; this is proved in the online notes of Milne~\cite{milne-modular-forms-notes} on modular forms and modular functions:
\begin{theorem}[Milne~\cite{milne-modular-forms-notes}, Theorem~6.1]
\label{milne-theorem}
The field~$\C(X_0(N))$ of modular functions for~$\Gamma_0(N)$ is generated over~$\C$ by~$j(q)$ and~$j(q^N)$. In particular, every modular function for~$\Gamma_0(N)$ is a rational function of~$j(q)$ and~$j(q^N)$.
\end{theorem}
By standard properties of modular forms, the quotient~$f/g$ of two modular forms of weight~$k$ for a congruence subgroup~$\Gamma$  is a modular function of weight~0 for that congruence subgroup, with possible poles at the zeroes of~$g$.

We therefore take~$f \in M_k(\Gamma_0(N))$ and consider the modular function
\begin{eqnarray}
\label{f-cdot-stuff}
g:=\frac{f \cdot E_4^a \cdot E_6^b}{\Delta^c} \in \M_0(\Gamma_0(N)),
\end{eqnarray}
where we choose~$a \in \{0,1,2\}$, $b \in \{0,1\}$ and~$k \in \N$ such that $k+4a+6b-12c=0$. This has weight~0 by definition, so by Theorem~\ref{milne-theorem} we can write it as a rational function of~$j(q)$ and~$j(q^N)$, which are both rational functions of~$\eta$-products. From before, we recall that~$E_4$, $E_6$ and~$\Delta$ are all rational functions of $\eta$-products, so we can write~$f$ as a rational function of $\eta$-products, which is what we set out to do.
\end{proof}
It is interesting to note that we only needed to use $\eta$-products of level~$4N$ here (the~$N$ coming from the level of our form~$f$, and the~4 coming from the $\eta$-quotient representations of~$E_4$ and~$E_6$). This corresponds to the figures in Table~\ref{table-of-spaces}, where we needed to raise the level by a factor of at least~4, corresponding to the level~4 $\eta$-quotients that we used for~$E_4$ and~$E_6$.

We see also that we have to have either the square of a prime or two distinct primes dividing the level, because we always need to write~$E_4$ and~$E_6$ in terms of $\eta$-products, and we need to find at least two different solutions to the set of simultaneous equations in the powers of the~$\eta^i$ (one for~$E_4$ and one for~$E_6$). This means that we need to have at least~3~$i$, so in particular we cannot write~$E_4$ and~$E_6$ just in terms of $\eta$-quotients of level~$p$ for any prime~$p$.

If we are willing to raise the level then we can use~$g:=\eta(q)^2\eta(q^{11})^2 \in S_2(\Gamma_0(11))$ as a multiplier in~\eqref{f-cdot-stuff}, so the only thing stopping our rational function being a polynomial is the numerator of the modular function of weight~0. However, we \emph{cannot} assume that~$g$ only has a zero at~$\infty$, because while the modular form~$\Delta(q) \in S_{12}(\sltz)$ only has a simple zero at~$\infty$ and no zeroes elsewhere, the modular form~$\Delta(q) \in S_{12}(\Gamma_0(N))$ (viewed as an oldform) has many zeroes; it has one at each cusp of~$\Gamma_0(N)$.

\section{Acknowledgements}
The author would like to thank Paul-Olivier Dehaye, Ken McMurdy and Alan Lauder for helpful conversations. The computations for this article were carried out on William Stein's computer {\tt MECCAH} and on the Mathematics Institute in Oxford's computers; the author is grateful for this computer time.

\end{document}